\documentclass[11pt]{article}
 
\usepackage{amssymb,epsfig} 
\usepackage{amsmath}
\usepackage{mathrsfs}

\renewcommand\mathcal{\mathscr}

\newtheorem{defi}{Definition}%[section]
\newtheorem{theo}[defi]{Theorem}
\newtheorem{coro}[defi]{Corollary}
\newtheorem{rema}[defi]{Remark}

\newcommand{\btheo}{\begin{theo}}
\newcommand{\etheo}{\end{theo}}
\newcommand{\brema}{\begin{rema}}
\newcommand{\erema}{\end{rema}}

\newcommand{\bcoro}{\begin{coro}}
\newcommand{\ecoro}{\end{coro}}
\newcommand{\dem}{\noindent{\bf Proof. }}
\newcommand{\rem}{\noindent{\bf Remark. }}

\newcommand{\R}{{\mathcal R}}
\newcommand{\E}{{\mathcal E}}
\renewcommand{\H}{{\mathcal H}}
\renewcommand{\O}{{\mathcal O}}
\newcommand{\I}{{\mathcal I}}
\renewcommand{\P}{{\mathcal P}}

\newcommand{\maths}[1]{{\mathbb #1}}

\newcommand{\OO}{\maths{O}}
\newcommand{\RR}{\maths{R}}
\newcommand{\NN}{\maths{N}}
\newcommand{\CC}{\maths{C}}
\newcommand{\QQ}{\maths{Q}}
\newcommand{\HH}{\maths{H}}
\newcommand{\KK}{\maths{K}}
\newcommand{\ZZ}{\maths{Z}}
\newcommand{\PP}{\maths{P}}

\newcommand{\ra}{\rightarrow}

\newcommand{\ov}[1]{{\overline #1}} 
\newcommand{\wt}[1]{{\widetilde{#1}}}
\newcommand{\wh}[1]{{\widehat{#1}}}

\newcommand{\ga}{\gamma}
\newcommand{\Ga}{\Gamma}

\newcommand{\cqfd}{\hfill$\Box$}

\title{On the closedness of approximation spectra}

\author{Jouni Parkkonen \and  Fr\'ed\'eric Paulin}

\renewcommand{\Re}{{\operatorname{Re}}}

\newcommand{\Heis}{{\operatorname{Heis}}}

\begin{document}

\maketitle 

 \begin{abstract}
Le spectre classique de Lagrange pour l'approximation des nombres
r\'eels par des rationnels, est ferm\'e, par un th\'eor\`eme de
Cusick. Plus g\'en\'eralement, nous montrons que de nombreux spectres
d'approximation sont ferm\'es, en utilisant des propri\'et\'es de
p\'en\'etration du flot g\'eod\'esique dans des voisinages de
pointes de vari\'et\'es \`a courbure strictement n\'egative, et un
r\'esultat de Maucourant  \cite{Mau}.

\bigskip

Generalizing Cusick's theorem on the closedness of the classical
Lagrange spectrum for the approximation of real numbers by rational
ones, we prove that 
various approximation spectra are closed, using
penetration properties of the geodesic flow in cusp neighbourhoods in
negatively curved manifolds and a result of Maucourant \cite{Mau}.
%\footnote{{\bf Keywords:} geodesic flow, negative curvature, Lagrange
%spectrum, Diophantine approximation, quaternions, Heisenberg
%group.~~{\bf AMS codes:} 53 C 22, 11 J 06}
\end{abstract}

\bigskip
\bigskip
\bigskip
The {\it approximation constant} of an irrational real number $x$ by
rational numbers is 
$$
c(x)=\liminf_{p,q\in\ZZ,\;q\ra+\infty}\; |q|^2\,\Big|x-\frac{p}{q}\Big|
$$ 
(though some references consider $c(x)^{-1}$ or even
$(2c(x))^{-1}$). The {\it Lagrange spectrum} $\operatorname{Sp}_\QQ$
is the subset of $\RR$ consisting of the $c(x)$ for
$x\in\RR-\QQ$. Many properties of $\operatorname{Sp}_\QQ$ are known
(see for instance \cite{CF}), and have been known for a very long
time, through the works of Korkine-Zolotareff, Hurwitz, Markoff, Hall,
\dots.  The fact that $\operatorname{Sp}_\QQ$ is a closed subset of
$\RR$ was proved by Cusick only relatively recently, in 1975.

For many examples of a locally compact ring $\wh K$ containing a dense
countable subring $K$, a linear algebraic group $\underline{G}$
defined over $K$ and a left invariant distance $d$ on the locally
compact group $\underline{G}(\wh K)$, one can define a similar
approximation spectrum of elements of $\underline{G}(\wh K)$ by
elements of $\underline{G}(K)$. In this note, we prove that
many such approximation spectra also are closed subsets of $\RR$, in
particular for 
\begin{itemize}
\item[$\bullet$] the approximation of complex numbers by elements in
  imaginary quadratic number fields,
\item[$\bullet$] the approximation of real Hamiltonian quaternions by
  rational ones, and
\item[$\bullet$] the approximation of elements of a real Heisenberg
  group by rational points.
\end{itemize}
In each of the above cases, the approximating elements are restricted to
certain subclasses of the quadratic irrational or rational elements, as
explained below. 
 
These arithmetic results will follow from a theorem in Riemannian
geometry, that we will state and prove, after recalling some
definitions.

\medskip
Let $M$ be a complete Riemannian manifold with dimension at least $2$
and sectional curvature at most $-1$, which is geometrically finite
(see for instance \cite{Bow} for a general reference).  Let $e$ be a
{\it cusp} of $M$, i.e.~an asymptotic class of minimizing geodesic
rays along which the injectivity radius goes to $0$.  In particular,
when $M$ has finite volume (which is going to be the case in all our
arithmetic applications), it is geometrically finite, and moreover,
the set of cusps of $M$ is in natural bijection with the (finite) set
of ends of $M$ (see loc.~cit.). Let $\rho_e:[0,+\infty[\;\ra M$ be a
    minimizing geodesic ray in $M$ in the class $e$ and let
    $\beta_e:M\ra \RR$ be Busemann's {\it height function} relative to
    $e$ (see for instance \cite[p.~268]{BH}) defined by
$$
\beta_e(x)=\lim_{t\ra+\infty} t-d(x,\rho_e(t))\;.
$$ 
Note that if another representative $\rho'_e$ of $e$ is considered,
then the new height function $\beta'_e$ only differs from $\beta_e$ by
an additive constant. 

Recall that a (locally) geodesic line $\ell:\RR\ra M$ {\it starts
  from} (resp.~{\it ends at}) $e$ if the map from $]a,+\infty[$ to
$M$, for some $a$ big enough, defined by $t\mapsto \ell(-t)$
(resp.~$t\mapsto \ell(t)$), is a minimizing geodesic ray in the class
$e$. A geodesic line $\ell$ is {\it positively recurrent} if there
exists a compact subset $K$ of $M$ and a sequence $(t_n)_{n\in\NN}$ in
$[0,+\infty[$ converging to $+\infty$ such that $\ell(t_n)\in K$ for
every $n$. For every positively recurrent geodesic line $\ell$
starting from $e$, define the {\it asymptotic height} of $\ell$ (with
respect to $e$) to be $\limsup_{t\ra+\infty} \beta_e(\ell(t))$.
Define (see for instance \cite{HPMZ,HPsurv}) the {\it asymptotic
  height spectrum} of $(M,e)$ as the set of asymptotic heights of
positively recurrent geodesic lines starting from $e$. If $C$ is a
compact subset of $M$, define the {\it height} of $C$ (with respect to
$e$) as
$$
\operatorname{ht}_e(C)=\max\{\beta_e(x)\;:\;x\in C\}\;.
$$ 
Note that the asymptotic height of a geodesic line, the asymptotic
height spectrum of $(M,e)$ and the height of a closed geodesic depend
on the choice of $\rho_e$ only up to a uniform additive constant.
There is a canonical normalization, by asking that $\rho_e(0)$ belongs
to the boundary of the maximal Margulis neighbourhood of $e$, see
\cite{BK,HPMZ} for instance. In some cases, however, this is not an
optimal choice in terms of computation lengths.

\medskip
Theorem \ref{theo:asympheispecclos} answers a question raised during
the work of the second author with S.~Hersonsky, see for instance page
233 in \cite{PPapp}.  In its proof, we will use the following result
of F.~Maucourant \cite[Theo.~2 (2)]{Mau}, whose main tool is Anosov's
closing lemma (and which builds on a partial result of \cite{HPMZ}). We
denote the unit tangent bundle of a Riemannian manifold $M$ by
$\pi:T^1M\to M$. A unit tangent vector is {\it periodic} if it is
tangent to a closed geodesic.

\btheo \label{theo:maucourant} Let $V$ be a complete Riemannian
manifold with sectional curvature at most $-1$, let
$(\phi^t)_{t\in\RR}$ be its geodesic flow, and let $J_0$ be the subset
of $T^1V$ of periodic unit tangent vectors. If $f:T^1V\ra \RR$ is a
proper continuous map, then
$$
\RR\cap\{\;\limsup_{t\ra+\infty} f(\phi^t v)\;:\;v\in T^1V\}=
\overline{\{\;\max_{t\in\RR} f(\phi^t v)\;:\;v\in J_0\}}\;.
\;\;\;\mbox{\cqfd}
$$  
\etheo

Here is the main geometric result of this note:

\btheo \label{theo:asympheispecclos}
The asymptotic height spectrum of $(M,e)$ is closed. It is equal to the 
closure of the set of the heights of the closed geodesics in $M$.  
\etheo

Notice that by \cite[Theo.~3.2 and 3.4]{HPMZ}, the asymptotic height
spectrum has a finite lower bound. But by \cite[Prop.~4.1]{HPMZ}, its
infimum (which is a minimum by Theorem \ref{theo:asympheispecclos}) is
not always attained by the height of a closed geodesic. Hence, its
minimum is not always isolated.

\medskip \dem Busemann's height function $\beta_e$ is continuous (in
fact $1$-Lipschitz). Let us prove that it is proper.  Let $\wt M\ra M$
be a universal Riemannian cover of $M$, with covering group $\Ga$. Let
$\partial_\infty \wt M$ be the sphere at infinity of $\wt M$, and
endow $\wt M\cup\partial_\infty\wt M$ with the cone topology (see for
instance \cite[p.~263]{BH}.  Let $\Lambda\Ga$ be the limit set of
$\Ga$ and let $\Omega \Ga=\partial_\infty\wt M-\Lambda\Ga$ be the
domain of discontinuity of $\Ga$ (see for instance \cite{Bow} and notice that
$\Omega\Gamma$ is empty if $M$ has finite volume): The group $\Ga$
acts properly discontinuously on $\wt M\cup\Omega \Ga$. The set of
ends of the quotient manifold with boundary $M^*=\Ga\backslash(\wt
M\cup\Omega \Ga)$ is in one-to-one bijection with the set of cusps of
$M$, by the map which to a minimizing geodesic ray defining a cusp of
$M$ associates the end of $M$ it converges to (see loc.~cit.). By
construction, $\beta_e(x)$ converges to $+\infty$ when $x$ converges
to the end of $M^*$ corresponding to $e$, and tends to $-\infty$ when
$x$ tends to any other end or any boundary point of $M^*$. This
implies that $\beta_e$ is proper.

Hence, the map $f=\beta_e\circ\pi:T^1M\ra \RR$ is also continuous and
proper. Note that if a geodesic line $\ell$ in $M$ is not positively
recurrent, then $\ell(t)$ converges, as $t$ goes to $+\infty$, either to
an end of $M^*$ or to a boundary point of $M^*$ (see loc.~cit.). Hence,
$\limsup_{t\ra+\infty} \beta_e(\ell(t))=\pm\infty$. Therefore, Theorem
\ref{theo:asympheispecclos} follows from Maucourant's Theorem
\ref{theo:maucourant}, applied to $f=\beta_e\circ\pi$. 
\cqfd

\bigskip
For our arithmetic applications, we transform Theorem
\ref{theo:asympheispecclos} into a form which is more applicable,
using the framework of Diophantine approximation in negatively curved
manifolds introduced in \cite{HPMZ,HPsurv}.  We recall the relevant
definitions from these references:

Let $\xi_e$ be the point at infinity of a lift $\wt\rho_e$ to $\wt M$
of the previously chosen minimizing geodesic ray $\rho_e$. Let
$\wt\beta_e:\wt M\ra\RR$ be Busemann's height function associated to
$\wt\rho_e$, i.e.~$\wt\beta_e(x)=\lim_{t\ra+\infty} t-
d(x,\wt\rho_e(t))$. A {\it horoball} (resp.~{\it horosphere}) centered
at $\xi_e$ is the preimage by $\wt\beta_e$ of $[s,+\infty[$
(resp. $\{s\}$) for some $s\in\RR$. A horoball $H$ centered at $\xi_e$
is {\it precisely invariant} under the action of the stabilizer
$\Ga_\infty$ of $\xi_e$ in $\Ga$ if the interiors of $H$ and $\ga H$
do not meet for any $\ga\in\Ga-\Ga_\infty$.  Let $H_e$ be a precisely
invariant horoball centered at $\xi_e$, which exists by \cite{Bow}.
Assume without loss of generality that $\wt\rho_e$ starts in $\partial
H_e$.

Let $\R_e$ be the set of geodesic lines in $M$ starting from and
ending at the cusp $e$, whose first point at height $0$ is at time
$0$. We endow $\R_e$ with its Fr\'echet filter. For every $r$ in
$\R_e$, define $D(r)$ as the length of the subsegment of $r$ between
its first and its last point whose height is $0$. By
\cite[Rem.~2.9]{HPMZ}, points $r\in\R_e$ go out of every finite subset
if and only if $D(r)$ tends to $+\infty$. Let $\operatorname{Lk}_e$ be
the set of positively recurrent geodesic lines starting from $e$,
whose first point at height $0$ is at time $0$.

For every distinct $x,y$ in $\operatorname{Lk}_e\cup\R_e$, define the
{\it cuspidal distance} $d'_e(x,y)$ between $x$ and $y$ as follows:
Let $\wt x$ be a lift of $x$ starting from $\xi_e$; for every $t>0$,
let $\H_t$ be the horosphere centered at $\wt x(+\infty)$, at signed
distance $-\log 2t$ from $\partial H_e$ along $\wt x$; then
$d'_e(x,y)$ is the minimum, over all lifts $\wt y$ of $y$ starting
from $\xi_e$, of the greatest lower bound of $t>0$ such that $\H_t$
meets $\wt y$ (see \cite[Sect. 2.1]{HPMZ}). The map $d'_e$ is an
actual distance in all our arithmetic applications (see loc.cit.), and
it depends on the choice of $\rho_e$ only up to a positive
multiplicative constant.

For every $x$ in $\operatorname{Lk}_e$, define the
{\it approximation constant} $c_{M,e}(x)$ of $x$ by elements of $\R_e$ as
$$
c_{M,e}(x)=\liminf_{r\in \R_e} \;\;d'_e(x,r)\;e^{D(r)}\;,
$$ and the {\it Lagrange spectrum} of $(M,e)$ as the subset of $\RR$
consisting of the constants $c_{M,e}(x)$ for $x$ in
$\operatorname{Lk}_e$.

\bcoro\label{coro:lagrangeferm}
The Lagrange spectrum of $(M,e)$ is closed.
\ecoro

\dem By \cite[Theo.~3.4]{HPMZ} (see also \cite[p.~232]{PPapp}), the
map $t\mapsto -\log(2t)$ is a homeomorphism from the Lagrange spectrum
onto the asymptotic height spectrum.  \cqfd

\medskip\noindent \rem
By the remark following Theorem \ref{theo:asympheispecclos}, the Lagrange
spectrum is bounded but its maximum is not always isolated.

\bigskip
Let us now give some arithmetic applications of this corollary, using
the notations introduced in \cite{PPCRAS,PPGT}.

Let $m$ be a squarefree positive integer, and let $\I$ be a nonzero
ideal of an order $\O$ in the ring of integers $\O_{-m}$ of the
imaginary quadratic number field $K_{-m}=\QQ(i\sqrt{m})$. For
$p_1,\dots,p_k\in\O$, let $\langle p_1,\dots,p_k\rangle$ be the ideal
of $\O$ generated by $p_1,\dots,p_k$.  Let 
$$
\E_{\I}= \{(p,q)\in\O\times\I \,:\; \langle p,q\rangle=\O\}\;.
$$ 
For every $x\in\CC-K_{-m}$, define the {\it approximation constant} of
$x$ by  elements of $\O\I^{-1}$ as
$$
c_\I(x)=\liminf_{(p,q)\,\in\,\E_{\I}\,,\;|q|\ra\infty} \;\; 
|q|^2\,\Big|x-\frac{p}{q}\Big|
$$
(the condition $\langle p,q\rangle=\O$ is not
needed when $\O$ is principal, for instance if $\O=\O_{-m}$ for $m=1,2, 3,
7, 11, 19, 43, 67,163$.) 
Define the {\it Bianchi-Lagrange spectrum} for the
approximation of complex numbers by elements of $\O\I^{-1}\subset
K_{-m}$ as the subset $\operatorname{Sp}_\I$ of $\RR$ consisting of
the $c_\I(x)$ for $x\in\CC-K_{-m}$.

\btheo \label{theo:corpquadimagi} 
The Bianchi-Lagrange spectrum $\operatorname{Sp}_\I$ is closed.  
\etheo

When $\I=\O=\O_{-m}$, this result is due to Maucourant \cite{Mau}.

\medskip
\dem Let $X=\HH^3_\RR$ be the upper halfspace model of the real
hyperbolic space of dimension $3$ (and sectional curvature $-1$).  The
group $\operatorname{SL}_2(\CC)$ acts isometrically on $X$, so that
its continuously extended action on $\partial_\infty X =
\CC\cup\{\infty\}$ is the action by homographies. Let $\Ga$ be the
(discrete) image in $\operatorname{Isom}(X)$ of the preimage of the
upper-triangular subgroup by the canonical morphism
$\operatorname{SL}_2(\O)\ra\operatorname{SL}_2(\O/\I)$. Let $\P_\Ga$
be the set of parabolic fixed points of elements of $\Ga$. Let
$M=\Ga\backslash X$, and let $e$ be its cusp corresponding to the
parabolic fixed point $\infty$ of $\Ga$. Note that $M$ is not
necessarily a manifold, as $\Ga$ may have torsion.  However, Theorem
\ref{theo:asympheispecclos} extends to this situation without any
changes.

By standard results in arithmetic subgroups (see for instance
\cite{BHC,Bor}, and the example (1) in \cite[\S 6.3]{PPGT}), $M$ has
finite volume and we have
$$
\P_\Ga=K_{-m}\cup\{\infty\}\,,
$$ so that $\partial_\infty X- \P_\Ga=\CC-K_{-m}$. Let $\Ga_\infty$ be
the stabilizer in $\Ga$ of the point $\infty$, which preserves the
Euclidean distance in $\partial_\infty X- \{\infty\}=\CC$. 

By \cite[Lem.~2.7]{HPMZ}, the map, which to $r\in \R_e$ associates the
double class modulo $\Ga_\infty$ of an element $\ga_r\in
\Ga-\Ga_\infty$ such that $\ga_r\infty$ is the other point at infinity
of a lift of $r$ to $X$ starting from $\infty$, is a bijection
$$
\R_e\ra\Ga_\infty\backslash (\Ga-\Ga_\infty)/\Ga_\infty\;.
$$ 
The map, which to $x\in\partial_\infty X- \{\infty\}$ associates
the image $\ell_x$ in $M$ of the geodesic line starting from $\infty$,
passing at Euclidean height $1$ at time $t=0$, and ending at $x$,
induces a bijection
$$ 
\Ga_\infty\backslash (\partial_\infty X-\P_\Ga)\ra
\operatorname{Lk}_e\;,
$$ 
since $\ell_x$ is positively recurrent if and only if $x\notin
\P_\Ga$ (see for instance \cite{Bow}). Furthermore (see for instance
\cite[page 314]{EGM}), the map, which to $(p,q)\in\E_\I$ associates
the image $r_{p/q}$ in $M$ of the geodesic line starting from
$\infty$, passing at Euclidean height $1$ at time $t=0$, and ending at
$p/q$, induces a bijection
$$ 
\Ga_\infty\backslash \Big\{\,\frac pq\;:\;(p,q)\in \E_\I\Big\} \ra \R_e\;.
$$ 
The horoball $\H_\infty$ of points with Euclidean height at least $1$
in $X$ is precisely invariant, since $\left(\begin{array}{cc}1 & 1\\
    0&0\end{array}\right)$ maps to an element of $\Ga$ and by Shimizu's
Lemma (see also \cite[\S 5]{HPMZ}). Let $\rho_e$ be the image by the
canonical projection $X\ra M$ of a geodesic ray from a point of
$\partial\H_\infty$ to $\infty$.  We use this minimizing geodesic ray
to define Busemann's height function $\beta_e$ and the cuspidal
distance $d'_e$. Hence, by definition, for every $r$ in $\R_e$, we
have
$$
D(r)=d_X(\H_\infty,\ga_r\H_\infty)\;.
$$ 
If $q$ is the lower-left entry of a representative in
$\operatorname{SL}_2(\CC)$ of an element $\ga$ in $\Ga-\Ga_\infty$,
then we have
$$
d_X(\H_\infty,\ga\H_\infty)=2\log |q|
$$ 
by \cite[Lem.~2.10]{HPMZ}. Hence, for every $(p,q)\in\E_\I$, we
have $D(r_{p/q})=2\log |q|$. 

It has been proved in \cite[\S 2.1]{HPMZ} (for the real hyperbolic
space $X$ of any dimension) that, for every $x,y$ in
$\operatorname{Lk}_e$, the cuspidal distance $d'_e(x,y)$ is equal to
the minimum of the Euclidean distances between the other points at
infinity of two lifts to $X$ of $x,y$ starting from $\infty$.

{}From the above, it follows that, for every $x\in \partial_\infty
X-\P_\Ga=\CC-K_{-m}$, 
$$
c_{M,e}(\ell_x)=\liminf_{(p,q)\in \E_\I,\,|q|\ra+\infty} 
e^{D(r_{p/q})}\;d'_e(\ell_x,r_{p/q})= c_\I(x)\;.
$$  
Hence, Theorem \ref{theo:corpquadimagi} follows from Corollary
\ref{coro:lagrangeferm}. \cqfd

\bigskip
Let $\I'$ be a nonzero two-sided ideal in an order $\O'$ of a
quaternion algebra $A(\QQ)$ over $\QQ$ ramifying over $\RR$, for
instance the Hurwitz ring $\O'=\ZZ[\frac{1}{2}(1+i+j+k),i,j,k]$ in
Hamilton's quaternion algebra over $\QQ$ with basis $(1,i,j,k)$, and
let $N$ be the reduced norm on $A(\RR)=A(\QQ)\otimes_\QQ\RR$ (see for
instance \cite{Vig}). Consider the set $$\E_{\I'}= \{(p,q)\,\in\,
\O'\times\I'\;:\; \exists\; r,s\,\in\,\O',\; N(qr-qpq^{-1}s)=1\}.$$
For every $x\in A(\RR)-A(\QQ)$, define the {\it approximation
constant} of $x$ by elements of $\O'\I'^{-1}\subset A(\QQ)$ as
$$
c_{\I'}(x)=\liminf_{(p,q)\,\in\,\E_{\I'}\,,
\;N(q)\ra\infty} \;\;  N(q)N(x-pq^{-1})^{\frac{1}{2}}\;,
$$ and the {\it Hamilton-Lagrange spectrum} for the approximation of
elements of $A(\RR)$ by elements of $\O'\I'^{-1}\subset A(\QQ)$ as the
subset $\operatorname{Sp}_{\I'}$ of $\RR$ consisting of the
$c_{\I'}(x)$ for $x\in A(\RR)-A(\QQ)$.

\medskip
\btheo \label{theo:quaternionentiers} 
The Hamilton-Lagrange spectrum $\operatorname{Sp}_{\I'}$ is closed.  
\etheo 

\dem The proof is the same as the previous one, with the following
changes.
\begin{itemize}
\item[$\bullet$] Let $X=\HH^5_\RR$ be the upper halfspace model of the
  real hyperbolic space of dimension $5$ (and sectional curvature
  $-1$).  With $\HH$ the field of quaternions of Hamilton, identified
  with $\RR^4$ by its standard basis $1,i,j,k$, we have
  $\partial_\infty X =\HH\cup\{\infty\}$. The group
  $\operatorname{SL}_2(\HH)$, of $2\times 2$ matrices with
  coefficients in $\HH$ and Dieudonn\'e determinant $1$, acts
  isometrically on $X$, so that its continuously extended action on
  $\partial_\infty X$ is, with the obvious particular cases,
  $$\big(\left(\begin{array}{cc}a& b \\ c &
      d\end{array}\right),z\big)\mapsto (az+b)(cz+d)^{-1},
  $$see for
  instance \cite{Kel}.
\item[$\bullet$] Let $\Ga$ be the image in $\operatorname{Isom}(X)$ of
  the preimage of the upper-triangular subgroup by the canonical
  morphism $\operatorname{SL}_2(\O')\ra\operatorname{SL}_2(\O'/\I')$.
\item[$\bullet$] Fix an identification of the quaternion algebras
  $A(\RR)$ and $\HH$. We have $\P_\Ga=A(\QQ)\cup\{\infty\}$ by the
  example (3) in \cite[\S 6.3]{PPGT}, so that $\partial_\infty X-
  \P_\Ga=\HH-A(\QQ)$.
\item[$\bullet$] The fact that the horoball $\H_\infty$ of points with
  Euclidean height at least $1$ in $X$ is precisely invariant is
  proved in \cite[page 1091]{Kel}.
\item[$\bullet$] By definition of $\E_{\I'}$ and of the Dieudonn\'e
  determinant, the map, which to $(p,q)\in\E_{\I'}$ associates the
  image $r_{pq^{-1}}$ in $M=\Ga\backslash X$ of the geodesic line
  starting from $\infty$, through $\partial\H_\infty$ at time $t=0$
  and ending at $pq^{-1}$, induces a bijection $ \Ga_\infty\backslash
  \{pq^{-1}\;:\;(p,q)\in\E_{\I'}\}\ra \R_e$.
\item[$\bullet$] If $q$ is the lower-left entry of a representative of
  an element $\ga$ in $\Ga-\Ga_\infty$, then we have
  $d_X(\H_\infty,\ga\H_\infty)=\log N(q) $ by \cite[Lem.~6.7]{PPGT},
  so that $D(r_{pq^{-1}})=\log N(q)$.
\item[$\bullet$] Recall that the reduced norm $N$ on $\HH$ is the
  square of the Euclidean distance on $\HH$ making the basis
  $(1,i,j,k)$ orthonormal. \cqfd
\end{itemize}

\bigskip Our last result concerns Diophantine approximation in
Heisenberg groups.  For every integer $n\geq 2$, consider the Lie
group
$$
\Heis_{2n-1}(\RR)=\{(z,w)\in\CC\times \CC^{n-1}\;:\; 2\; \Re\; z
-|w|^2=0\}\;,
$$ 
where $w'\cdot\overline{w}=\sum_{i=1}^{n-1} w_i'\ov {w_i}$ is the
standard Hermitian scalar product on $\CC^{n-1}$ and
$|w|^2=w\cdot\overline{w}$, with the group law 
$$
(z,w)(z',w')= (z + z' + w'\cdot \overline{w}, w + w')\;.
$$ 
Consider the {\it modified Cygan distance} $d'_{\rm Cyg}$ on
$\Heis_{2n-1}(\RR)$, which is defined (uniquely)  as the distance which
is invariant
under left translations and satisfies
$$
d'_{\rm Cyg}((z,w),(0,0))= \sqrt{2\,|z|+|w|^2}\;,
$$ see \cite[\S 6.1]{PPGT}. Notice that its induced length distance is
equivalent to the Cygan distance and to the Carnot-Carath\'eodory
distance (see \cite{Gol}).

Let $\I$ be a nonzero ideal of an order $\O$ in the ring of
integers $\O_{-m}$ of the imaginary quadratic number field
$K_{-m}=\QQ(i\sqrt{m})$, and let $\omega$ be an element of $\O_{-m}$
with $\operatorname{Im} \omega>0$ such that $\O=\ZZ+\omega \ZZ$. Notice
that $\Heis_{2n-1}(\RR)$ is the set of real points of a $\QQ$-form
$\Heis_{2n-1}$ (depending on $m$) of the $(2n-1)$-dimensional Heisenberg
group, whose set of $\QQ$-points is $\Heis_{2n-1}(\RR)\cap
(K_{-m}\times {K_{-m}}^{n-1})$.

If $n=2$ and $\O=\O_{-m}$, then let $\E'_{\I}$ be the set of $(a,
\alpha, c)\in \O_{-m}\times \I\times \I$ such that $2\;{\rm Re}\; a
\overline{c} = |\alpha|^2$ and $\langle a,\alpha,c\rangle =\O_{-m}$.
Otherwise, see the fifth point below for the definition of
$\E'_{\I}$. For every $x\in\Heis_{2n-1}(\RR)-\Heis_{2n-1}(\QQ)$,
define the {\it approximation constant} $c'_{\I}(x)$ of $x$ by
$$
c'_{\I}(x)=\liminf_{(a,\alpha,c)\in\,\E'_{\I}\,,\; |c|\ra\infty} \;\; 
|c|\;d'_{\rm Cyg}(x,(a/c,\alpha/c))\;,
$$ 
and the {\it Heisenberg-Lagrange spectrum} for the approximation of
elements of $\Heis_{2n-1}(\RR)$ by elements of $\{(a/c,\alpha/c)\;:\;
(a,\alpha,c)\in\,\E'_{\I}\}\subset \Heis_{2n-1}(\QQ)$ as the subset of
$\RR$ consisting of the $c'_{\I}(x)$ for $x\in\Heis_{2n-1}(\RR) -
\Heis_{2n-1}(\QQ)$.

\medskip
\btheo \label{theo:heisenberg} 
The Heisenberg-Lagrange spectrum is closed.  
\etheo

\dem The proof is the same as that of Theorem
\ref{theo:corpquadimagi}, with the following changes.
\begin{itemize}
\item[$\bullet$] Let $X=\HH^n_\CC$ be the Siegel domain model of the
  complex hyperbolic $n$-space, which is the manifold
  $$
  \{(w_0,w)\in\CC\times\CC^{n-1}\;:\;2\operatorname{Re}
  w_0-|w|^2>0\}
  $$ 
  with the Riemannian metric
$$ 
ds^2=\frac{
(dw_0-dw\cdot\overline{w})(\overline{dw_0}-
w\cdot\overline{dw})+(2\operatorname{Re}
w_0-|w|^2)\;dw\cdot\overline{dw}}
{(2\operatorname{Re}
  w_0-|w|^2)^2}
$$ 
(we normalized the metric so that the maximal sectional curvature
is $-1$). Its boundary at infinity is $\partial_\infty X=
\Heis_{2n-1}(\RR)\cup\{\infty\}$.
\item[$\bullet$] Using matrices by blocks in the decomposition
  $\CC^{n+1}= \CC\times\CC^{n-1}\times\CC$ with coordinates
  $(z_0,z,z_n)$, let $Q$ be the matrix of the Hermitian form
  $-z_0\overline{z_n} - z_n\overline{z_0}+|z|^2$ of signature $(1,n)$,
  and let $\operatorname{SU}_Q$ be the group of complex matrices of
  determinant $1$ preserving this Hermitian form. We identify
  $X\cup\partial_\infty X$ with its image in the complex projective
  $n$-space $\PP_n(\CC)$ by the map (using homogeneous coordinates)
  $(w_0,w)\mapsto [w_0:w:1]$ and $\infty\mapsto[1:0:0]$. The group
  $\operatorname{SU}_Q$, acting projectively on $\PP_n(\CC)$, then
  preserves $X$ and acts isometrically on it.
\item[$\bullet$] Let $\Ga$ be the image in $\operatorname{Isom}(X)$ of
  the preimage, by the canonical morphism from
  $\operatorname{SU}_Q\cap \operatorname{SL}_{n+1}(\O)$ to
  $\operatorname{SL}_{n+1}(\O/\I)$, of the subgroup of matrices
  all of whose coefficients in the first column vanish except the
  first one.
\item[$\bullet$] We have $\P_\Ga=\Heis_{2n-1}(\QQ)\cup\{\infty\}$ by
  the example (2) in \cite[\S 6.3]{PPGT}, so that $\partial_\infty X-
  \P_\Ga=\Heis_{2n-1}(\RR)-\Heis_{2n-1}(\QQ)$.
\item[$\bullet$] Let us consider the horoball
  $$ 
  \H_\infty=\{(w_0,w)\in\CC\times\CC^{n-1}\;:\;2\operatorname{Re}
  w_0-|w|^2\geq 4\operatorname{Im} \omega\}\;.
  $$  
  The fact that $\H_\infty$ is precisely invariant is proved in
  \cite[Lem.~6.4]{PPGT}.
\item[$\bullet$] Let $\E'_{\I}$ be the set of
  $(a,\alpha,c)\in\O\times{\I}^{n-1}\times\I$ such that there exists a
  matrix of the form $\left(\begin{array}{ccc}a & \ga & b\\ \alpha & A
  & \beta \\ c&\delta &d \end{array}\right)$ that defines an element
  of $\Ga$. If $n=2$ and $\O=\O_{-m}$, we recover the previous
  notation, by \cite[\S 6.1]{PPGT}. By definition, the map, which to
  $(a,\alpha,c)\in\E'_{\I}$ associates the image in $M$ of the
  geodesic line starting from $\infty$, through $\partial\H_\infty$ at
  time $t=0$ and ending at $(ac^{-1},\alpha c^{-1})$, induces a
  bijection
  $$
  \Ga_\infty\backslash \big\{(ac^{-1},\alpha c^{-1})\;:\;
  (a,\alpha,c)\in\E'_{\I}\big\}\ra \R_e.
  $$
\item[$\bullet$] If $c$ is the lower-left entry of an element $\ga$ in
  $\Ga-\Ga_\infty$, then we have $d_X(\H_\infty,\ga\H_\infty)=\log
  |c|+\log(2\operatorname{Im} \omega) $ by \cite[Lem.~6.3]{PPGT}.
\item[$\bullet$] By \cite[Prop.~6.2]{PPGT}, the cuspidal distance is
  equal to a multiple of the modified Cygan distance. Hence, there
  exists a constant $\kappa>0$ such that $c'_{\I}(x)=\kappa
  \,c_{M,e}(\ell_x)$ for every $x\in \partial_\infty X-\P_\Ga$.
  \cqfd
\end{itemize}

\medskip
Other applications could be obtained by varying the nonuniform
arithmetic lattices in $\operatorname{Isom}(\HH^n_\KK)$ with
$\KK=\RR,\CC,\HH,\OO$ (where $n=2$ in this last octonion case).

\medskip
%{\small 
{\it Acknowledgements}: We thank S.~Hersonsky for 
conversations during the April 2008 Workshop on Ergodic Theory and Geometry
at the University of Manchester.
%}

\noindent {\small
\begin{tabular}{l}
Department of Mathematics and Statistics, P.O. Box 35\\
40014 University of Jyv\"askyl\"a, FINLAND\\
{\it e-mail: parkkone@maths.jyu.fi}
\end{tabular}
\\
\smallskip
 \mbox{}
\\
\begin{tabular}{l}
D\'epartement de Math\'ematique et Applications, UMR 8553 CNRS\\
\'Ecole Normale Sup\'erieure, 45 rue d'Ulm\\
75230 PARIS Cedex 05, FRANCE\\
{\it e-mail: Frederic.Paulin@ens.fr}
\end{tabular}
}

\end{document}